\documentclass[a4paper,11pt]{article}

\usepackage{a4wide}
\usepackage{amsthm,amsmath,amssymb}

\def\C{\mathbb C}
\def\cond{{\rm cond}}
\def\cdt{\!\cdot\!}
\def\eps{\varepsilon}
\def\gam{\gamma}
\def\Hom{{\rm Hom}}
\def\GL{{\rm GL}}
\def\Ind{{\rm Ind}}
\def\ext{{\rm \bf ext}}
\def\N{\mathbb N}
\def\OF{{\mathcal O}}
\def\Q{\mathbb Q}
\def\st{\rm St}
\def\res{{\rm \bf res}}

\newtheorem{theo}{Theorem}
\newtheorem{lemma}{Lemma}[section]
\theoremstyle{definition}
\newtheorem{rem}[lemma]{Remark}

\begin{document}

\title{Test vectors for trilinear forms : the case of two principal series}
\author{Mladen Dimitrov and Louise Nyssen \cr {\footnotesize
 dimitrov@math.jussieu.fr,  lnyssen@math.univ-montp2.fr} } 
\date{\today}

\maketitle

\section{Introduction}

\bigskip
Let $F$ be a finite extension of  $\Q_p$ with ring of integers $\OF$
and uniformizing parameter $\pi$. Let  $V_1$, $V_2$ and
$V_3$ be three irreducible, admissible, infinite dimensional  
representations of  $G=\GL_2(F)$ of central characters 
$\omega_1$, $\omega_2$ and $\omega_3$ and conductors $n_1$, $n_2$ and $n_3$. 
Using the theory of Gelfand pairs, Diprenda Prasad proves in \cite{P}
that the space of $G$-invariant linear forms on $V_1\otimes V_2
\otimes V_3$  has dimension at most one and gives 
a precise criterion for this dimension to be  one, that we will now explain.

Let $D^*$ be the group of invertible elements of the unique quaternion
division algebra  $D$ over $F$.  
When $V_i$ is a discrete series representation of $G$, denote
by $V'_i$  the irreducible representation of $D^*$
associated  
to $V_i$ by the Jacquet-Langlands correspondence. Again, by
the theory of Gelfand pairs,  
the space of $D^*$-invariant linear forms on $V'_1\otimes V'_2 \otimes
V'_3$ has dimension at most one. 

A necessary condition for the existence on a  non-zero 
$G$-invariant linear form on $V_1\otimes V_2\otimes V_3$
(resp. non-zero $D^*$-invariant linear form on $V'_1\otimes V'_2 \otimes
V'_3$), that we will {\it always assume}, is that 
$$ \enspace \omega_1\omega_2\omega_3=1.$$

Let $\sigma_i$ be the two dimensional representations of the
Weil-Deligne group of $F$ associated to  $V_i$.  
The triple tensor product $\sigma_1 \otimes \sigma_2 \otimes\sigma_3$
is an eight dimensional symplectic representation of the Weil-Deligne
group having a local root number $\eps(\sigma_1 \otimes \sigma_2
\otimes \sigma_3)$ equal to  $1$ or $-1$. 
When $\eps(\sigma_1 \otimes \sigma_2 \otimes \sigma_3)=- 1$, 
one can prove that the  $V_i$'s are all discrete series representations of $G$.

\begin{theo} (Prasad  \cite[Theorem 1.4]{P})  
If all the  $V_i$'s are supercuspidal, 
assume that the residue characteristic of $F$ is not 2. Then \par
$\centerdot$ $\eps(\sigma_1 \otimes \sigma_2 \otimes \sigma_3)=1$ if,
and only if, there exists a non-zero $G$-invariant linear form on  
$V_1\otimes V_2 \otimes V_3$ , and \par 
$\centerdot$ $\eps(\sigma_1 \otimes \sigma_2 \otimes \sigma_3)=-1$ if,
and only if,  there exists a non-zero $D^*$ invariant linear form on  
$V'_1\otimes V'_2 \otimes V'_3$.\par
\end{theo}

Given a non zero $G$-invariant linear form $\ell$ on $V_1\otimes V_2
\otimes V_3$, or a non-zero  $D^*$-invariant linear form $\ell'$ on  
$V'_1\otimes V'_2 \otimes V'_3$, the goal is to find a vector in  $V_1
\otimes V_2 \otimes V_3$ which is not in the kernel of $\ell$,  
or a vector in  $V'_1\otimes V'_2 \otimes V'_3$ which is not in the
kernel of $\ell'$.  Such a vector is called a test vector.  
The following results of Prasad and Gross-Prasad show that new vectors 
can sometimes be used as test vectors. In what follows 
$v_i$ denotes a new vector in $V_i$ (see \S\ref{nv}). 

\begin{theo}\label{vt-000} (Prasad \cite[Theorem 1.3]{P}) 
If all the $V_i$'s are unramified principal series, then  $v_1 \otimes v_2 \otimes v_3$ is a test vector.  
\end{theo} 

\begin{theo}\label{vt-111} (Gross and Prasad  \cite[Proposition  6.3]{GP}) 
Suppose  all   the $V_i$'s are  unramified twists of the Steinberg
representation. 
\begin{itemize}
\item If $\eps(\sigma_1 \otimes \sigma_2 \otimes \sigma_3)=1$,
then $v_1 \otimes v_2 \otimes v_3$ is a test vector. 
\item If $\eps(\sigma_1 \otimes \sigma_2 \otimes \sigma_3)=-1$
and if $R$ is the the unique maximal order in $D$, then any vector
belonging to  the unique line in  $V'_1\otimes V'_2 \otimes V'_3$
fixed by $R^* \times R^*\times R^*$ is a test vector.
\end{itemize}
\end{theo} 

Actually, the proof by  Gross and Prasad of the first statement of the
above theorem contains another result : 

\begin{theo}\label{vt-110} 
If two of the $V_i$'s are  unramified twists of the Steinberg 
representation and the third 
one is an  unramified principal series, then $v_1 \otimes v_2
\otimes v_3$ is a test vector.  
\end{theo}

However, as mentioned in \cite{GP}, new vectors are not always
test vectors. 
Let $K=\GL (\OF)$ be the maximal compact subgroup of $G$ and suppose
that $V_1$ and $V_2$ are unramified, but  $V_3$ is ramified. 
Since $v_1$ and $v_2$ are $K$-invariant and $\ell$ is $G$-equivariant,
$v \mapsto  \ell(v_1\otimes v_2\otimes v)$ defines a $K$-invariant
linear form on $V_3$. Since $V_3$ is ramified, so is its contragredient, 
and  therefore the above linear form has to vanish. In particular
$\ell(v_1\otimes v_2\otimes v_3)=0.$ 

To go around this obstruction for new vectors to be test vectors, Gross
and Prasad made the following suggestion : suppose that $V_3$ has
conductor $n=n_3\geq 1$; since $V_3$ has unramified central character,
its contragredient  
representation has non-zero invariant vectors by the $n$-th standard Iwahori subgroup 
$I_n=\begin{pmatrix} \OF^\times & \OF \\ \varpi^n\OF &  \OF^\times\end{pmatrix}$ of 
$G$; put  $\gam=\begin{pmatrix}\pi^{-1}  & 0
  \\ 0 &  1\end{pmatrix}$ and let  $v_1^*\in V_1$ be a non-zero vector
on the line  fixed by
the maximal compact subgroup $\gam^{n}K\gam^{-n}$ of $G$; 
since $K\cap \gam^{n}K\gam^{-n}=I_n$, the  linear form on $V_3$
given by $v \mapsto  \ell(v_1^*\otimes v_2\otimes v)$ is not
necessarily zero  and there is still hope for $v_1^* \otimes v_2
\otimes v_3$ to be   a test vector. This is the object of the 
following theorem

\begin{theo}\label{vt-00n} If $V_1$ and $V_2$  are unramified and 
$V_3$ has conductor $n_3$, then $v_1^* \otimes v_2
\otimes v_3$ is a test vector, where $v_1^*=\gam^{n_3}\cdt v_1$. 
\end{theo}
Theorem  \ref{vt-00n} for  $n_3=1$, together with Theorems \ref{vt-000},
\ref{vt-111} and \ref{vt-110},  
completes the study of test vectors when the $V_i$'s have
conductors $0$ or $1$ and unramified central characters. 

\medskip
Assume from now on that $V_1$ and $V_2$ are (ramified or
unramified) principal series. 
Then for  $i=1,2$ there exist  quasi-characters  $\mu_i$ and $\mu'_i$
of  $F^\times$ such that $\mu'_i\mu_i^{-1}\neq |\cdot|^{\pm 1}$, and
 $$V_i = \Ind_{B}^{G} \chi_i \text{ ,  with }
 \chi_i  \begin{pmatrix} a & b \cr 0 & d \cr \end{pmatrix}
  =  \mu_i(a)\mu'_i(d).$$
  
According to Theorem 1 there exists a non-zero $G$-invariant linear
form $\ell$ on $V_1\otimes V_2 \otimes V_3$, so 
we are looking for a test vector in $V_1\otimes V_2 \otimes V_3$.
 The following theorem  is  our main result.

\begin{theo}\label{vt-mkn}
Suppose that $V_1$ and $V_2$  are principal series
such that  $\mu_1$ and $\mu'_2$ are unramified. Put 
$$x=\max(n_2-n_1,n_3-n_1) \qquad{\rm and}\qquad v_1^*=\gam^{x}\cdt v_1.$$
Then $x\geq 0$  and, if  $v_1^* \otimes v_2 \otimes v_3$ is {\it not} a test vector, then
\begin{itemize}
\item either $n_1=0$, $n_2=n_3>0$ and  $\gam^{n_2-1}\cdt v_1
\otimes v_2 \otimes v_3$ is  a test vector, 
\item or $n_2=0$, $n_1=n_3>0$ and  
$v_1 \otimes \gam\cdt v_2 \otimes v_3$ is  a test vector,
\item or $\widetilde{V_3}$ is a quotient of $\Ind_{B}^{G}(\chi_1 \chi_2 \delta^{\frac{1}{2}})$, 
$n_1+n_2=n_3$ and $ v_1 \otimes \gam^{n_1}\cdt v_2 \otimes v_3$ is  a test vector. 
\end{itemize} 
\end{theo}
The assumptions of the theorem imply in particular
that  $V_1$ and $V_2$ have minimal conductor among their twists.
If $V_1$ and $V_2$  are two arbitrary principal series, then 
one can always find characters  $\eta_1$, $\eta_2$ and $\eta_3$ of $F^\times$
with $\eta_1 \eta_2 \eta_3 =1$,  such that the above theorem applies to 
$(V_1 \otimes \eta_1) \otimes(V_2\otimes \eta_2 ) \otimes(V_3\otimes
\eta_3)$. Nevertheless, we found also interesting to study the case
when $\mu_1$ or  $\mu'_2$ is ramified. Then we are able to show that 
certain new vectors are {\it not} test vectors, while {\it a priori} 
this cannot be seen by a direct argument
(the  obstruction 
of  Gross and Prasad  described above does not apply to this case). 
Put $m_1=\cond(\mu'_1)$ and $m_2=\cond(\mu'_2)$

\begin{theo}\label{no-vt} Suppose that  $\mu_1$ or  $\mu'_2$ is   ramified. Let  $x$, $y$ and $z$ be integers such that 
\begin{itemize}
\item $x\geq m_1$,
\item $y \geq m_2$,
\item $x-n_3\geq z \geq y$, and 
\item  $x-y\geq \max( n_1-m_1, n_2-m_2,1)$. 
\end{itemize}
Put 
\begin{equation}\label{v12}
\begin{split}
v_1^*=\begin{cases} \gam^{x-m_1}\cdt v_1 &  \text{ , if }  \mu'_1 \text{ is ramified,}\\
\gam^{x}\cdt v_1-\beta_1 \gam^{x-1}\cdt v_1  & \text{ , if }  \mu'_1 \text{ is unramified.}
\end{cases}  \\
v_2^*=\begin{cases} \gam^{y-m_2}\cdt v_2 &  \text{ , if }  \mu_2 \text{ is ramified.}\\
\gam^{y-n_2}\cdt v_2-\alpha_2^{-1} \gam^{y-n_2+1}\cdt v_2  & \text{ , if }  \mu_2 \text{ is unramified.}
\end{cases}
\end{split}
\end{equation}
Then 
$$\ell(v_1^* \otimes  v_2^* \otimes \gam^{z} \cdt v_3)=0.$$
\end{theo}

\bigskip
We will prove theorems \ref{vt-mkn} and \ref{no-vt} by following  the 
pattern of the proof of Theorem \ref{vt-000} in \cite{P}, with the necessary changes. 

%

\bigskip
We believe that suitable generalization of the 
method of Gross and Prasad would give  test vectors in the 
case where at least two of the  $V_i$'s are 
special representations, as well as in the case where
one is a special representation and one is a principal series. 
 On the other hand in order to find test vectors
in the case where  at least two of the  $V_i$'s are supercuspidal,
one should use different techniques, involving probably computations in 
Kirillov models.

\bigskip
The search for test vectors in our setting is motivated by subconvexity problems for 
$L$-functions of triple products of automorphic forms on $\GL(2)$. 
Roughly speaking, one wants to bound the value of the  $L$-function along the
critical line $\Re(z)=\frac{1}{2}$. 
In \cite{BR1} and \cite{BR2} Joseph Bernstein and Andre Reznikov  establish a  
{\it subconvexity bound} when the {\it eigenvalue} attached to one of the representations varies. 
Philippe Michel and Akshay Venkatesh 
considered the case when the {\it level} of one representation varies. 
More details about subconvexity and those related techniques can be
found in \cite{V} or \cite{MV}.  
Test vectors are key ingredients. Bernstein and Reznikov use an explicit test vector. 
Venkatesh uses a theoretical one, but explains that the bounds would
be better with an explicit one (see \cite[\S 5]{V}).  

There is an extension of Prasad's result in \cite{HS}, where Harris and Scholl prove that 
the dimension of the space of $G$-invariant linear forms on $V_1\otimes V_2 \otimes V_3$ is one when $V_1$, $V_2$ and $V_3$ 
are principal series representations (either irreducible or reducible, but with infinite dimensional irreducible subspace).  
They apply their result to the global setting to construct elements in the motivic cohomology of the product of two modular curves predicted by Beilinson.

\subsection*{Acknowledgments.}
We would like to  thank Philippe Michel for suggesting the study of this problem, 
and of course Benedict Gross and Diprenda Prasad for their articles full of inspiration.
The second named author would like to  thank also Paul Broussous and Nicolas Templier for many interesting discussions, and  Wen-Ching Winnie Li for the opportunity to spend one semester at PennState University where the first draft of this paper was written. 
  
\section{Background on induced admissible representations of $\GL(2)$.}

\subsection{About induced and contragredient representations.}\label{notations}
Let $(\rho, W)$ be a smooth representation of a closed subgroup $H$ of
$G$. Let $\Delta_H$ be the modular function on $H$.  
The induction of $\rho$ from $H$ to $G$,  denoted $\Ind_{H}^{G} \rho
$,  is the space   of functions $f$ from $G$ to $W$ 
satisfying the two following conditions :

(1) $\forall h \in H, \quad \forall g \in G, \quad
f(hg)=\Delta_H(h)^{-\frac{1}{2}} \rho(h) f(g)$, 

(2) there exists an open compact subgroup $K_f$ of $G$ such that
$$\forall k \in K_f, \quad \forall g \in G, \quad f(gk)= f(g)$$ 

\noindent where  $G$ acts by right translation as follows : 
$$\forall g, g' \in G, (g\cdot f)(g') = f(g'g).$$ 
With the  additional condition  that $f$ must be compactly supported modulo $H$, 
one gets the {\it compact} induction denoted by ${\rm ind}_{H}^{G}$. 
When $G/H$ is compact, there is no difference between 
$\Ind_{H}^{G}$ and ${\rm ind}_{H}^{G}$.  

Let $B$ the Borel subgroup of upper triangular matrices in $G$, and let $T$ be the diagonal torus. 
The character  $\Delta_T$ is trivial and we will use $\Delta_B=\delta^{-1} $ with 
$\delta \begin{pmatrix} a & b \cr 0 & d \cr \end{pmatrix}  = \vert
\frac{a}{d} \vert$  
where $ \vert \enspace \vert$ is the normalised valuation of $F$. The
quotient $B \backslash G$ is compact and can be identified with $\mathbb{P}^1(F)$.

For a smooth representation $V$ of $G$, the contragredient
representation $\widetilde{V}$ is the space of smooth linear forms
$l$ on $V$, where $G$  acts as follows : 
$$\forall g\in G, \qquad  \forall v\in V, \qquad (g\cdot l)(v)= l(g^{-1}\cdot v).$$ 

We refer the reader to  \cite{BZ} for more details about induced and contragredient representations.

\bigskip
\subsection{New vectors and ramification.}\label{nv}
Let $V$ be an irreducible, admissible, infinite dimensional
representation of  $G$ with central character $\omega$.  Then
$\widetilde{V}\cong V\otimes \omega^{-1}$.
To  the descending chain of  compact subgroups of $G$ 
$$ K=I_0 \supset I_1 \supset \cdots \supset I_n \supset I_{n+1} \cdots  $$
one can associate an ascending chain of vector spaces
$$V^{I_0,\omega}=V^K \text{ , and for all }  n\geq 1, \quad 
V^{I_{n},\omega}=\left\{v\in V \Big{|}
\begin{pmatrix} a & b \\ c &
      d\end{pmatrix}\cdt v=\omega(d)v \text{ , for all }
\begin{pmatrix} a & b \\ c & d\end{pmatrix}\in I_n \right\}.$$
There exists a minimal $n$ such that the vector space
$V^{I_{n},\omega}$ is non-zero.  
It is necessarily one dimensional and any non-zero vector in it is called a {\it new
  vector} of $V$.    
The integer $n$ is the {\it conductor} of $V$. 
The representation $V$ is said to be {\it unramified} if $n=0$. 

More information about new vectors can be found in \cite{C}.

\bigskip

\subsection{New vectors as functions on $G$.}\label{nv-functions}

Let $V$ be a principal series of  $G$, with central character
$\omega$, and conductor $n$.  
There exist  quasi-characters  $\mu$ and $\mu'$
of  $F^\times$ such that $\mu'\mu^{-1}\neq |\cdot|^{\pm 1}$, and
 $$V = \Ind_{B}^{G} (\chi)  \qquad \text{with } \qquad
 \chi  \begin{pmatrix} a &  * \cr 0 & d \cr \end{pmatrix}
  =  \mu(a)\mu'(d).$$
  
Then $\omega=\mu\mu'$ and $n = \cond(\mu)+\cond(\mu')$. A new vector $v$ in $V$ is  a non-zero function from $G$ to $\C$ such that   
for all $b \in B$, $g \in G$ and $k= \begin{pmatrix} * &  * \\ * & d  \end{pmatrix} \in I_n$ 
$$v(bgk)=\chi(b)\delta(b)^\frac{1}{2} \omega(d)v(g).$$   
Put $$\alpha^{-1}= \mu(\pi)|\pi|^{\frac{1}{2}} \qquad{\rm and} \qquad \beta^{-1}= \mu'(\pi)|\pi|^{-\frac{1}{2}}.$$

First, we assume that $V$ is unramified, and we normalise $v$ so that $v(1)=1$. 

\begin{lemma}\label{calcul-NR} If $V$ is unramified then for all $r\in \N$, 
$$(\gam^{r}\cdt v)(k)=
\begin{cases} \beta^r  & \text{ , if } k\in K\backslash I, \\
\alpha^s\beta^{r-s}  & \text{ , if } k\in I_s\backslash I_{s+1}
\text{  for  } 1\leq s\leq r-1,\\
\alpha^{r} & \text{ , if } k\in I_r.\\
\end{cases}$$

Similarly,  
$$(\gam^{r}\cdt v-\alpha^{-1}\gam^{r+1}\cdt v)(k)=
\begin{cases} \alpha^s\beta^{r-s}-\alpha^{s-1}\beta^{r+1-s} & 
\text{ , if } k\in I_s\backslash I_{s+1}
\text{  for  } 0\leq s\leq r,\\
0 & \text{ , if } k\in I_{r+1}.\\
\end{cases} $$

Finally, for $r\geq 1$,
$$(\gam^{r}\cdt v-\beta\gam^{r-1}\cdt v)(k)=
\begin{cases} \alpha^r(1-\frac{\beta}{\alpha}) & 
\text{ , if } k\in I_r,\\
0 & \text{ , if } k\in K \backslash I_{r}.
\end{cases} $$
\end{lemma}

\noindent{\it Proof : }
If $k\in I_r$, then $\gam^{-r}k\gam^{r}\in K$, so
$$(\gam^{r}\cdt v)(k)=\alpha^{r}v(\gam^{-r}k\gam^{r})= \alpha^{r}.$$ 

Suppose that $k =\begin{pmatrix}a & b \\ c & d\end{pmatrix}\in
I_s\backslash I_{s+1}$ for some $0\leq s\leq r-1$ (recall that $I_0=K$).
Then $\pi^{-s} c \in {\OF}^{\times}$ and
$$ (\gam^{r}\cdt v)(k)=\alpha^{r}v\begin{pmatrix}a & \pi^{r}b \\ \pi^{-r}c & d\end{pmatrix} 
= \alpha^{r}v\begin{pmatrix} (ad-bc)\pi^{r-s} & a \\ 0 & \pi^{-r}c \end{pmatrix}
= \alpha^s\beta^{r-s}.$$ 
The second part of the lemma follows by a direct computation. $\hfill \Box $ 

\medskip

For the rest of this section we assume that $V$ is ramified, that is $n \geq 1$. 
We put $$m=\cond(\mu') \qquad {\rm so \quad that} \qquad n-m=\cond(\mu).$$ 

By Casselman \cite[pp.305-306]{C} the restriction to $K$ of a new vector 
$v$ is supported by the double coset of 
$\begin{pmatrix} 1 & 0 \\ \pi^{m} & 1  \end{pmatrix}$ modulo $I_{n}$. 
In particular if  $\mu'$ is unramified ($m=0$), then $v$ is supported by 
$$I_{n} \begin{pmatrix} 1 & 0 \\ 1 & 1  \end{pmatrix}I_{n}
=I_{n} \begin{pmatrix} 0 & 1 \\ 1 & 0  \end{pmatrix}I_{n} =K\backslash I.$$
If $1 \leq m \leq n-1$, then $v$ is supported by 
$$I_{n} \begin{pmatrix} 1 & 0 \\ \pi^m & 1  \end{pmatrix}I_{n}
=I_{m} \backslash I_{m+1}.$$
If  $\mu$  is unramified ($m=n$), then $v$ is supported by  $I_{n}$. 
We normalise $v$ so that
$$v\begin{pmatrix} 1 & 0 \\ \pi^{m} & 1  \end{pmatrix}=1.$$

\begin{lemma}\label{calcul-SP0} 
If $\mu$ and $\mu'$ are both ramified ($0<m<n$), 
then for all $r\in \N$ and $k\in K$, 
$$(\gam^{r}\cdt v)(k)=
\begin{cases}\alpha^r \mu\Bigl(\frac{\det{k}}{\pi^{-(m+r)}c}\Bigr)\mu'(d) &
  \text{ , if } 
k=\begin{pmatrix} * & * \\ c & d  \end{pmatrix}\in I_{m+r}\backslash I_{m+r+1}, \\
0  & \text{ ,  otherwise}.
\end{cases}$$
\end{lemma}

\noindent{\it Proof : } 
For $k=\begin{pmatrix} a & b \\ c & d
\end{pmatrix}\in K$ we have 
$$\alpha^{-r}(\gam^{r}\cdt v)(k)= v(\gam^{-r}k\gam^{r})= v\begin{pmatrix} a &\pi^r b \\ \pi^{-r}c & d  \end{pmatrix}.$$

It is easy to check that for every $s\geq 1$, 
$$K\cap B\gam^{r}I_s\gam^{-r}=I_{s+r}.$$ 
It follows that $\gam^{r}\cdt v$ has its support in $I_{m+r}\backslash I_{m+r+1}$.
If $k\in I_{m+r}\backslash I_{m+r+1}$ then  $c \in {\pi}^{m+r} {\OF}^{\times}$, $d\in   {\OF}^{\times}$ and  we have the following decomposition :
\begin{equation}\label{decomposition}
\begin{pmatrix} a &\pi^r b \\ \pi^{-r}c & d  \end{pmatrix}=
\begin{pmatrix} \det{k} & \pi^{-m}cb\\ 0 & \pi^{-m-r}cd  \end{pmatrix}
\begin{pmatrix} 1 & 0 \\ \pi^{m} & 1  \end{pmatrix}
\begin{pmatrix} d^{-1} & 0 \\ 0 & \pi^{m+r}c^{-1}  \end{pmatrix}.
\end{equation}

Hence $$\alpha^{-r}(\gam^{r}\cdt v)(k)=
\mu\Bigl(\det(k)\Bigr)\mu'(\pi^{-m-r}cd)(\mu\mu')( \pi^{m+r}c^{-1})=
\mu\Bigl(\frac{\det(k)}{\pi^{-(m+r)}c}\Bigr)\mu'(d).$$ 
$\hfill \Box $

Similarly we obtain : 

\begin{lemma}\label{calcul-SP1} 
Suppose that $\mu$ is  unramified and $\mu'$ is ramified. Then, for all $r\in \N$ and $k\in K$, 
$$ (\gam^{r}\cdt v)(k)=
\begin{cases}\alpha^r \mu'(d) & \text{ , if } 
k=\begin{pmatrix} * & * \\ * & d  \end{pmatrix}\in I_{{n}+r}, \\
0  & \text{ ,  otherwise}.
\end{cases} $$
$$\Bigl( \gam^{r}\cdt v-\alpha^{-1} \gam^{r+1}\cdt v \Bigr)(k)=
\begin{cases}\alpha^r \mu'(d) & \text{ , if } 
k=\begin{pmatrix} * & * \\ * & d  \end{pmatrix}\in I_{n+r}\backslash I_{n+r+1}, \\
0  & \text{ ,  otherwise}.
\end{cases} $$
\end{lemma}

\begin{lemma}\label{calcul-SP2} 
Suppose that $\mu'$ is  unramified and $\mu$ is ramified. 
Then for all $r\in \N$, 
$$(\gam^{r}\cdt v)(k)=
\begin{cases}\alpha^s\beta^{r-s}
  \mu\left(\frac{\det(k)}{\pi^{-s}c}\right) & \text{ , if }  
k=\begin{pmatrix} * & * \\ c & *  \end{pmatrix}\in I_s\backslash
I_{s+1}\text{, with } 0\leq s\leq r,  \\
0  & \text{ , if }   k\in I_{r+1}.
\end{cases} $$

Moreover, if  $r \geq 1$, then  
$$\Bigl(\gam^{r}\cdt v-\beta \gam^{r-1}\cdt v \Bigr)(k) =
\begin{cases}\alpha^r
  \mu\left(\frac{\det(k)}{\pi^{-r}c}\right) & \text{ , if } k=\begin{pmatrix} * & * \\ c & *  \end{pmatrix}\in I_r\backslash I_{r+1},  \\
0  & \text{ ,  otherwise}.
\end{cases}
$$
\end{lemma}

\noindent{\it Proof : } We follow the pattern of proof of lemma \ref{calcul-SP0}.
The restriction of $\gam^{r}\cdt v$ to $K$ is zero outside
$$K\cap B\gam^{r}(K\backslash I)\gam^{-r}=K\backslash I_{r+1}.$$ 
For $ 0\leq s\leq r$ and 
$k=\begin{pmatrix} a & b \\ c & d  \end{pmatrix}\in I_s\backslash I_{s+1}$ 
we use the following decomposition :
\begin{equation}\label{decomposition2}
\begin{pmatrix} a &\pi^r b \\ \pi^{-r}c & d  \end{pmatrix}=
\begin{pmatrix} -\frac{\det{k}}{\pi^{-r}c} & a+\frac{\det{k}}{\pi^{-r}c}\\ 0 & \pi^{-r}c 
\end{pmatrix}
\begin{pmatrix} 1 & 0 \\ 1 & 1  \end{pmatrix}
\begin{pmatrix} 1 & 1+\frac{d}{\pi^{-r}c} \\ 0 & -1  \end{pmatrix}.
\end{equation}
Since $d\in\OF$ and $\pi^{r}c^{-1}\in\OF$ we deduce that : 
$$\alpha^{-r}(\gam^{r}\cdt v)(k)=
\mu\Bigl(\frac{\det{k}}{\pi^{-r}c}\Bigr)
\mu'(-\pi^{-r}c)\left\vert \pi^{r}c^{-1} \right\vert =
\mu\Bigl(\frac{\det{k}}{\pi^{-s}c}\Bigr)\alpha^{s-r}\beta^{r-s}.$$ 
$\hfill \Box $

As direct consequence of these lemmas we obtain 

\begin{lemma}\label{support} 
Let $v_1^*$ and $v_2^*$ be as in Theorem \ref{no-vt}. Then the support of $v_1^*$ is 
$$\begin{cases} 
I_{x}\backslash I_{x+1} &  
\text{ , if } \mu_1  \text{ is ramified, }\\
I_{x} &  
\text{ , if } \mu_1 \text{ is unramified, }
\end{cases}$$
and the support of $v_2^*$ is 
$$\begin{cases} 
I_{y}\backslash I_{y+1} &  
\text{ , if }    \mu'_2 \text{ is ramified, }\\
K \backslash I_{y+1} &  
\text{ , if }   \mu'_2 \text{ is unramified. }
\end{cases}$$
\end{lemma}

\section{Going down Prasad's exact sequence.}\label{suites-exactes}

In this section we will  explain how Prasad finds a non-zero 
$\ell\in \Hom_G ( V_1 \otimes
V_2 \otimes V_3 , \C)$ in the case
 where $V_1$ and $V_2$ are principal series representations.

\subsection{Prasad's exact sequence.}
The space  $\Hom_G ( V_1 \otimes V_2 \otimes V_3 , \C)$   
is canonically isomorphic to $\Hom_G (
V_1 \otimes V_2 , \widetilde{V_3} )$, hence  finding $\ell$ it is the
same as finding  a non-zero element $\Psi$ in it.  We have
$$V_1 \otimes V_2 = {\rm Res}_{G}\,\Ind_{B \times B}^{G \times G} \Bigl( \chi_1 \times \chi_2 \Bigr)$$
where the restriction is taken with respect to the 
diagonal embedding of $G$ in $G\times G$.  The action of $G$ on
$(B\times B) \backslash (G\times G) \cong \mathbb{P}^1(F)\times \mathbb{P}^1(F)$ has precisely two
orbits. 

The first is the diagonal $\Delta_{B \backslash G}$, which is closed and
can be identified with $B \backslash G$.  The second   is  its
complement   which is open and  can be identified with $T \backslash
G$ via the bijection : 
 $$ \begin{matrix}
T \backslash G & \longrightarrow & \Bigl( B \backslash G \times B \backslash G \Bigr) \setminus \Delta_{B \backslash G} \cr
\hfill Tg & \longmapsto & \Bigl( Bg , B\begin{pmatrix} 0 & 1 \cr 1 & 0
  \cr \end{pmatrix}g \Bigr) \hfill 
\end{matrix}$$

Hence, there is  a short exact sequence of $G$-modules :
\begin{equation}\label{courtesuite}
0 \rightarrow {\rm ind}_{T}^{G}\Bigl( \chi_1\chi'_2 \Bigr)  \xrightarrow{\ext} V_1 \otimes V_2 \xrightarrow{\res} 
\Ind_{B}^{G} \Bigl(\chi_1 \chi_2 \delta^{\frac{1}{2}}  \Bigr)\rightarrow 0,
\end{equation}
where $ \chi'_2  \begin{pmatrix} a & b \cr 0 & d \cr \end{pmatrix}
  =  \mu'_2(a)\mu_2(d).$
The surjection  $\res$ is given by the restriction to the diagonal. 
The injection $\ext$ takes a function 
$f \in {\rm ind}_{T}^{G}\Bigl( \chi_1\chi'_2 \Bigr)$ to a function 
$F \in \Ind_{B \times B}^{G \times G} \Bigl( \chi_1 \times \chi_2
\Bigr)$ vanishing on $\Delta_{B \backslash G}$, such that for all $g\in G$ 
$$F \Bigl( g, \begin{pmatrix} 0 & 1 \cr 1 & 0 \cr \end{pmatrix} g
\Bigr) = f(g) \label{rel}.$$

 Applying the functor $\Hom_G \Bigl( \bullet , \widetilde{V_3}
 \Bigr) $ yields  a long exact sequence :

\begin{multline}\label{longuesuite} 
0 \rightarrow \Hom_G \Bigl( \Ind_{B}^{G} \Bigl(\chi_1 \chi_2 \delta^{\frac{1}{2}}  \Bigr) , \widetilde{V_3} \Bigr) 
  \rightarrow  \Hom_G \Bigl( V_1 \otimes V_2 , \widetilde{V_3} \Bigr) 
  \rightarrow  \Hom_G \Bigl( {\rm ind}_{T}^{G}\Bigl( \chi_1\chi'_2 \Bigr), \widetilde{V_3} \Bigr) \\
\hfill \downarrow \hskip20mm\\
\hfill \cdots  \leftarrow {\rm Ext}_G^1 \Bigl( \Ind_{B}^{G} \Bigl(\chi_1 \chi_2 \delta^{\frac{1}{2}}  \Bigr) , \widetilde{V_3} \Bigr)
\end{multline}

\bigskip

\subsection{The simple case.}\label{cas-simple}

 The situation is easier if $V_3$  occurs in
$\Ind_{B}^{G} (\chi_1^{-1} \chi_2^{-1} \delta^{-\frac{1}{2}})$.
Then $\chi_1\chi_2$ does not factor through the determinant and there is a natural surjection  
$$\Ind_{B}^{G} \Bigl(\chi_1 \chi_2 \delta^{\frac{1}{2}}  \Bigr)
\twoheadrightarrow \widetilde{V_3}.$$
This surjection is an isomorphism, unless there exists a quasi-character
$\eta$ of $F^{\times}$ such that  $\chi_1\chi_2\delta=\eta \circ \det$ 
in which case the kernel is a line generated by the function $\eta \circ {\rm det}$. 
From (\ref{courtesuite}) we obtain 
a surjective  homomorphism $\Psi$ completing the following commutative
diagram :  
\begin{equation}\label{diagramme}
\begin{matrix}
 V_1 \otimes V_2  & \xrightarrow{\res}  & \Ind_{B}^{G} \Bigl(\chi_1 \chi_2 \delta^{\frac{1}{2}}  \Bigr) \cr
\hfill {\scriptstyle \Psi}\! \searrow  & & \swarrow \hfill\cr
&\widetilde{V_3}&\cr
\end{matrix}
\end{equation}

Finding a test vector is then reduced to finding an element of $V_1
\otimes V_2$ whose image by $\res$ is  not zero  (resp. not a multiple
of $\eta \circ {\rm det}$), if  
$V_3$ is  principal series (resp.  special representation). 

\medskip

Following the notations of paragraph \ref{nv-functions} put, for $i=1$ and $i=2$
$$m_i=\cond(\mu'_i) \qquad \alpha_i^{-1}= \mu_i(\pi)|\pi|^{\frac{1}{2}}\quad \text{and} \qquad 
 \beta_i^{-1}= \mu'_i(\pi)|\pi|^{-\frac{1}{2}}.$$

\subsubsection{Proof of theorem \ref{no-vt} in the simple case.}

To prove theorem \ref{no-vt}, suppose that  $\mu_1$ or   $\mu'_2$ is  ramified.
By our assumptions $x>y$, hence $I_x\cap (K\backslash I_{y+1})=\varnothing$. 
Therefore the supports of $v_1^*$ and $v_2^*$ are disjoint and 
$$\res( v_1^* \otimes v_2^*)=0.$$ 

Using the  diagram (\ref{diagramme}) we see that for any $v \in V_3$ : 
$$\ell( v_1^* \otimes v_2^* \otimes v ) 
= \Psi(v_1^* \otimes v_2^*)(v)=0.$$

In particular $\ell( v_1^* \otimes v_2^* \otimes \gam^z\cdt v_3)=0$ which 
proves Theorem \ref{no-vt} in the simple case. 

\bigskip
The rest of section \ref{cas-simple}  will be devoted to the proof of Theorems  \ref{vt-00n} 
and \ref{vt-mkn} in the simple case. 
Consequently, we will suppose that $\mu_1$ and $\mu'_2$ are unramified, that is $m_1-n_1=m_2=0$.

\subsubsection{Proof of Theorem \ref{vt-00n} in the simple case.}

Since $V_1$ and  $V_2$ are unramified, by theorem \ref{vt-000} we may
assume that  $V_3$ is ramified. Then necessarily 
$$\widetilde{V_3}= \eta \otimes \st, $$ 
where $\st$ is the Steinberg representation and $\eta$ is an
unramified character. Hence $ n_3=1$ 
and we will prove that $\gam \cdot v_1\otimes v_2\otimes v_3$ is a test vector. 
\bigskip

The function 
$$\left\{ \begin{matrix}
G & \longrightarrow & \C \hfill\\
g & \mapsto & \eta \Bigl({\rm det}(g)\Bigr)^{-1}\res(\gam \cdot v_1\otimes v_2)(g) \\
\end{matrix}\right.$$ 
is  not constant, since   according to lemma \ref{calcul-NR}  
$$  \eta \Bigl({\rm det}(1)\Bigr)^{-1}(\gam \cdot v_1 \otimes v_2)(1)  
= v_1( \gam)  v_2(1)=\alpha_1 $$
and
$$ \eta \Bigl({\rm det}\begin{pmatrix} 0 & 1 \\ 1 & 0  \end{pmatrix}\Bigr)^{-1}  (\gam \cdot v_1 \otimes v_2 )\begin{pmatrix} 0 & 1 \\ 1 & 0  \end{pmatrix}  
=  \eta(-1) v_1\begin{pmatrix} 1 & 0 \\ 0 & \pi^{-1} \end{pmatrix}
=  \beta_1, $$
and $\alpha_1\neq \beta_1$ because $V_1$ is a principal series. 

\medskip
Hence $\Psi(\gam \cdot v_1\otimes v_2)\neq 0$. 
Moreover, since $$\gam K \gam^{-1} \cap K =I$$ we deduce that 
$$\Psi(\gam \cdot v_1\otimes v_2)\in \widetilde{V_3}^{I,{\omega_3}^{-1}}.$$
Hence $\Psi(\gam \cdot v_1\otimes v_2)$ cannot vanish on the line
${V_3}^{I,{\omega_3}}$, which is generated by $v_3$, and therefore 
 $\gam \cdot v_1\otimes v_2 \otimes v_3$ is a test vector.
 
\medskip
This  completes the proof of Theorem \ref{vt-00n} in the simple case.

\subsubsection{Proof of Theorem \ref{vt-mkn} in the simple case, when
  $\widetilde{V_3}$  is  a special representation.} 

Assume now that  
$$\widetilde{V_3}= \eta \otimes \st, $$ 
where $\st$ is the Steinberg representation and $\eta$ is a character.
Since $$\eta = \mu_1\mu_2|\cdot|=\mu'_1\mu'_2|\cdot|^{-1}$$
and $\mu_1$ and $\mu'_2$ are unramified, it follows that 
$\eta$ is unramified if, and only if, both  $V_1$ and
$V_2$ are unramified. Since this case was taken care of in the
previous paragraph, we can assume for the rest of this paragraph that 
$\eta$ is  ramified. Then 
$$n_1=n_2=\cond(\eta)\geq 1 \qquad \text{and} \qquad n_3=2n_1=n_1+n_2.$$ 
We will prove that $v_1\otimes \gam^{{n_1}} \cdot v_2\otimes v_3$ is a test vector. 
\bigskip

The function 
$$\left\{ \begin{matrix}
G & \longrightarrow & \C \hfill\\
g & \mapsto & \eta \Bigl({\rm det}(g)\Bigr)^{-1}\res(v_1\otimes \gam^{{n_1}} \cdot v_2)(g) \\
\end{matrix}\right.$$ 
is  not constant, since according to lemmas \ref{calcul-SP1}  and \ref{calcul-SP2}  
$$  \eta \Bigl({\rm det}(1)\Bigr)^{-1}(v_1\otimes \gam^{{n_1}} \cdot v_2)(1)  = 0 $$
whereas
$$ \eta \Bigl({\rm det}\begin{pmatrix} 1 & 0 \\ \pi^{n_1} & 1  \end{pmatrix}\Bigr)^{-1}  
(v_1\otimes \gam^{{n_1}} \cdot v_2 )\begin{pmatrix} 1 & 0 \\ \pi^{n_1}
  & 1  \end{pmatrix}   
=  \alpha_2^{n_1} \neq 0. $$

\medskip
Hence $\Psi(v_1\otimes \gam^{{n_1}} \cdot v_2)\neq 0$. 
Moreover, since 
$$I_{n_1}\cap \gam^{n_1}I_{n_2}\gam^{-n_1}=I_{n_1+n_2}=I_{n_3}$$  
we deduce that 
$$\Psi(v_1\otimes \gam^{{n_1}} \cdot v_2)\in \widetilde{V_3}^{I_{n_3},{\omega_3}^{-1}}.$$
Hence  $\Psi(v_1\otimes \gam^{{n_1}} \cdot v_2)$ cannot vanish on the line
${V_3}^{I_{n_3},{\omega_3}}$, which is generated by $v_3$, and therefore 
  $v_1\otimes \gam^{n_1} \cdot v_2\otimes v_3$ is a test vector.

\subsubsection{Proof of Theorem \ref{vt-mkn} in the simple case, when
  $\widetilde{V_3}$ is a principal series.} 
\label{cas-simple2} 

Finally,  we  consider  the case  where $\widetilde{V_3}$ is a principal series representation. Then
$$ \widetilde{V_3} = \Ind_{B}^{G} \Bigl(\chi_1 \chi_2 \delta^{\frac{1}{2}}  \Bigr)$$ and 
$$ n_3 = \cond(\mu_1\mu_2)+ \cond(\mu'_1\mu'_2)=n_2+n_1.$$
We will prove that $v_1\otimes \gam^{n_1} \cdot v_2\otimes v_3$ is a test vector. 
\medskip

 According to lemmas \ref{calcul-NR}, 
\ref{calcul-SP1}  and \ref{calcul-SP2} we have  
$$ (v_1\otimes \gam^{n_1} \cdot v_2) \begin{pmatrix}1 & 0 \cr
  \pi^{n_1} & 1\cr\end{pmatrix} = \alpha_2^{n_1} \neq 0,$$
hence  $\res(v_1\otimes \gam^{n_1} \cdot v_2)\neq 0$. 

Therefore  $\Psi(v_1\otimes \gam^{n_1}v_2)\neq 0$.
Moreover, since 
$$I_{n_1}\cap \gam^{n_1}I_{n_2}\gam^{-n_1}=I_{n_1+n_2}=I_{n_3}$$  
we deduce that 
$$\Psi(v_1\otimes \gam^{n_1}v_2)\in (\widetilde{V_3})^{I_{n_3},{\omega_3}^{-1}}.$$
Hence $\Psi(v_1\otimes \gam^{n_1}v_2)$ 
cannot vanish on the line  ${V_3}^{I_{n_3},{\omega_3}}$, 
which is generated by $v_3$. 
Thus $v_1\otimes \gam^{n_1} \cdot v_2\otimes v_3$ is a test vector. 

This completes the proof of Theorem \ref{vt-mkn} in the simple case.

\subsection{The other case.}
The situation is more complicated if $\Hom_G ( \Ind_{B}^{G} (\chi_1
\chi_2 \delta^{\frac{1}{2}}  ) , \widetilde{V_3} ) =0$. By
 \cite[Corollary 5.9]{P} we have 
${\rm Ext}_G^1 ( \Ind_{B}^{G}(\chi_1 \chi_2 \delta^{\frac{1}{2}}
) , \widetilde{V_3} )=0$, hence  
 the long exact sequence (\ref{longuesuite}) yields the following  isomorphism :
$$\Hom_G \Bigl( V_1 \otimes V_2 , \widetilde{V_3} \Bigr) \simeq \Hom_G
\Bigl( {\rm ind}_{T}^{G}( \chi_1\chi'_2 ), \widetilde{V_3}
\Bigr).$$ 
Finally, by  Frobenius reciprocity   
$$\Hom_G \Bigl( {\rm ind}_{T}^{G}( \chi_1\chi'_2 ) , \widetilde{V_3} \Bigr) \simeq 
\Hom_T \Bigl( \chi_1\chi'_2  , \widetilde{V_{3\vert T}} \Bigr).$$
By \cite[ Lemmes 8-9]{W}  the latter space is one dimensional, 
since the restriction of $\chi_1\chi'_2$ to the center equals  
$\omega_3^{-1}$ (recall that $\omega_1\omega_2\omega_3=1$). 
Thus, we have four canonically isomorphic lines with corresponding
bases : 

\begin{equation}\label{chaine}
\begin{matrix}
0\neq\ell & \in & {\rm Hom}_G  \Bigl( V_1 \otimes V_2 \otimes V_3 , \C \Bigr) \\
             &     & \downarrow \wr \\       
0\neq\Psi     & \in & {\rm Hom}_G \Bigl( V_1 \otimes V_2 ,
\widetilde{V_3} \Bigr) \\     &     & \downarrow \wr \\    
0\neq\Phi    & \in & {\rm Hom}_G \Bigl( {\rm ind}_{T}^{G}(\chi_1\chi'_2), \widetilde{V_3} \Bigr) \\          &     & \downarrow \wr \\    
0\neq\varphi  & \in & {\rm Hom}_T \Bigl( \chi_1\chi'_2  , \widetilde{V_{3\vert T}} \Bigr) \\
\end{matrix}
\end{equation}
Observe that $\varphi$ can be seen as a linear form on $V_3$ satisfying  :
\begin{equation}\label{phi}
\forall t \in T, \qquad \forall v \in V_3, \qquad
\varphi(t\cdt v) = 
(\chi_1\chi'_2)(t)^{-1}\varphi(v).
\end{equation}

\begin{lemma}\label{lemmeV3}
$\varphi(v_3) \neq 0$ if, and only if, $\mu_1\mu'_2$ is unramified.
\end{lemma}

\noindent{\it Proof : } Suppose $\varphi(v_3) \neq 0$. Since $v_3\in V_3$ is a
new vector, for all $a,d\in\OF^\times$ we have 
$$\begin{pmatrix} a & 0 \\ 0 &
  d\end{pmatrix}\cdt v_3=\omega_3(d)v_3=
(\mu_1\mu'_1\mu_2\mu'_2)(d)^{-1}v_3.$$
Comparing it with (\ref{phi}) forces $\mu_1\mu'_2$ to be
unramified.

Conversely, assume that $\mu_1\mu'_2$ is unramified.
Take any $v\in V_3$ such
that $\varphi(v) \neq 0$. By smoothness $v$ is fixed by the 
principal congruence subgroup $\ker(K\rightarrow
\GL_2(\OF/\pi^s))$, for some $s\geq 0$. Then 
$\varphi(\gam^{s}\cdt v)=(\mu_1\mu'_2)(\pi^s)\varphi(v) \neq
0$ and $\gam^{s}\cdt v$ is fixed by the congruence subgroup 
$$I_{2s}^{1}:=\left\{k\in K \Big{|} k\equiv \begin{pmatrix} 1 & * \\ 0 &
      1\end{pmatrix} \pmod{\pi^{2s}} \right\}.$$
By replacing $\gam^{s}\cdt v $ by $v$ and $2s$ by $s$, we may assume
that $v\in V_3^{I_{s}^{1}}$ for some  $s\geq 0$. Since $I_{s}/I_{s}^{1}$
is a finite abelian group, $V_3^{I_{s}^{1}}$ decomposes as a direct sum 
of spaces 
indexed by the characters of  $I_{s}/I_{s}^{1}$. Then  $\varphi$ has
to be non-zero on 
$V_3^{I_{s},\omega_3}$ (defined in paragraph \ref{nv})
since by (\ref{phi}), $\varphi$ vanishes on all other summands of
$V_3^{I_{s}^{1}}$.  

By Casselman \cite[Theorem 1]{C} the space $V_3^{I_{s},\omega_3}$ has dimension $n_3-s+1$
and has a basis 
$$\Bigl( \quad v_3\quad, \quad \gam\cdt v_3\quad,  \dots , \quad\gam^{n_3-s}\cdt v_3 \quad\Bigr)$$  
(recall that $n_3$ denotes the conductor of $V_3$). 
Again by  (\ref{phi}),  $\varphi(\gam^{i}\cdt v_3)\neq 0$ for some
$i$ is equivalent to  $\varphi(v_3)\neq 0$. $\hfill \Box$

Notice that, when $\mu_1\mu'_2$ and  $\mu'_1\mu_2$ are both 
unramified, the claim follows from the first case in \cite[Proposition
2.6]{GP} applied  to the split torus $T$ of $G$.

\section{Going up Prasad's exact sequence.}

In this section we take as a starting point  lemma \ref{lemmeV3} 
and follow the isomorphisms (\ref{chaine}).

\subsection{From $\varphi$ to $\Phi$.}\label{phi-f}

Let  $x$, $y$ and $z$ be integers  such that 
$$x-n_3\geq z \geq y \geq 0 \qquad {\rm and} \qquad x-y \geq 1.$$ 
For the proof of Theorem \ref{vt-mkn} we will take  
$$x=\max(n_1,n_3)\geq 1 \qquad {\rm and} \qquad y=z=0.$$
Given a quasi-character $\mu$ of $F^\times$ define : 
$$\OF^{\mu}=\begin{cases}\OF & \text{ , if } \mu \text{ is unramified, }   \\
\OF^\times   & \text{ , if } \mu \text{ is ramified.} \end{cases}$$

Put $$ I_{f}=\begin{pmatrix} 1 & \pi^{-y}\OF^{\mu'_2} \\ 
 \pi^{x}\OF^{\mu_1}& 1\end{pmatrix}, $$ 
and consider the unique function  $f\in {\rm ind}_{T}^{G}( \chi_1\chi'_2 )$
which is zero outside the open compact subset $T I_f$ of $T\backslash
G$ and such that for all $b_0\in \pi^{-y}\OF^{\mu'_2}$
and $ c_0 \in \pi^{x}\OF^{\mu_1}$ we have : 
\begin{equation}\label{f}
f\begin{pmatrix} 1 & b_0 \\ c_0 &  1\end{pmatrix} = 
\begin{cases} 
\mu_1(\frac{\pi^{x}}{c_0})\mu'_2(b_0\pi^{y}) & \text{ , if } \mu_1
\text{ and }   \mu'_2 \text{ are ramified ; } \\

\mu'_2(b_0\pi^{y}) & \text{ , if }  \mu_1
\text{ is unramified } \text{ and }   \mu'_2 \text{ is ramified ;  } \\

\mu_1(\frac{\pi^{x}}{c_0}) & \text{ , if }\mu_1
\text{ is ramified } \text{ and }   \mu'_2 \text{ is unramified ; } \\ 

1 & \text{ , if } \mu_1
\text{ and }   \mu'_2 \text{ are unramified. }  \\
\end{cases}
\end{equation}

Since $x-n_3\geq z \geq y \geq 0$ and $x-y \geq 1$ we have
$$ I_{f}\subset \gam^{z}I_{n_3}^1\gam^{-z}$$
and so every   $k_0\in I_{f}$ fixes $\gam^{z}\cdt v_3$.

\medskip
By definition, the function   $g  \mapsto  f(g)\varphi( g\gam^{z}\cdt v_3)$
on $G$  factors through $T\backslash G$ and

$$\Bigl( \Phi(f) \Bigr)(\gam^{z}\cdt v_3) = 
\int_{T \backslash G} \! f(g) \, \varphi( g\gam^{z}\cdt v_3) dg
=\varphi( \gam^{z}\cdt v_3 ) \int_{ I_{f}} \! f(k_0) dk_0.$$ 

If we write $k_0=\begin{pmatrix} 1 & b_0 \\ c_0 & 1\end{pmatrix}\in
I_f$, then by separating the variables $b_0$ and  $c_0$ we obtain 

$$ \int_{ I_{f}} \!
f(k_0) dk_0=\begin{cases}  |\pi|^{x-y} & \text{, if } \mu_1 \text{ and }   \mu'_2 \text{ are unramified,} \\
0 & \text{, otherwise.}
\end{cases}$$
 
From this and from lemma \ref{lemmeV3} we deduce : 

\begin{lemma}\label{phi-f-v_3}
$\Phi(f)  (\gam^{z}\cdt v_3)\neq 0$ if, and only if, 
 $\mu_1$ and  $\mu'_2$ are both unramified.
 \end{lemma}

\subsection{From $\Phi$ to $\Psi$.} 

Now, we are going to compute $F=\ext(f)$ as a function on $G\times G$.
Recall that $F:G\times G\rightarrow \C$ is a function such that : 

- for all $b_1,b_2\in B$, $g_1,g_2\in G$,
$F(b_1g_1,b_2g_2)=\chi_1(b_1)\chi_2(b_2)\delta^{\frac{1}{2}}(b_1b_2)F(g_1,g_2)$, 

- for all $g\in G$,  $F(g,g)=0$ and $F(g,\begin{pmatrix} 0 & 1 \\ 1 &
  0\end{pmatrix}g) = f(g)$. 
    
\noindent Since $G=BK$, $F$ is uniquely determined by its restriction to
$K\times K$. 
Following the notations of paragraph \ref{nv-functions} put 
$$\alpha_i^{-1}= \mu_i(\pi)|\pi|^{\frac{1}{2}}\quad \text{and} \qquad 
 \beta_i^{-1}= \mu'_i(\pi)|\pi|^{-\frac{1}{2}}.$$ 

\begin{lemma} \label{FV} 
Suppose that  $x-n_3\geq z \geq y \geq 0$ and $x-y \geq\max (n_1-m_1, n_2-m_2,1)$. 
Then for  all $k_1 =\begin{pmatrix} * & * \\ c_1 & d_2 \end{pmatrix}$ 
and $k_2 =\begin{pmatrix} * & * \\ c_2 & d_2 \end{pmatrix}$ in $K$ we
have $F(k_1,k_2)=0$ unless 
$$d_1c_2 \neq 0, \qquad \frac{c_1}{d_1}\in \pi^x{\OF}^{\mu_1} \qquad {\rm and}\qquad  \frac{d_2}{c_2}\in \pi^{-y}{\OF}^{\mu'_2},$$ 
in which case, if we denote by
$s$ the valuation of $c_2$, we have

$$F(k_1,k_2)=\begin{cases}  
\mu_1\left(\frac{\det(k_1)}{\pi^{-x}c_1}\right) \mu'_1(d_1)
\mu_2\left(\frac{-\det(k_2)}{\pi^{-s}c_2}\right) \mu'_2(d_2)
\left(\frac{\alpha_2}{\beta_2}\right)^s 
& \text{, if }\mu_1 \text{ and }   \mu'_2 \text{ are ramified ; } \\ 
\mu'_1(d_1)
\mu_2\left(\frac{-\det(k_2)}{\pi^{-s}c_2}\right) \mu'_2(d_2)
\left(\frac{\alpha_2}{\beta_2}\right)^s
& \text{, if } \mu_1
\text{ is unramified } \text{ and }   \mu'_2 \text{ is ramified ;  } \\
\mu_1\left(\frac{\det(k_1)}{\pi^{-x}c_1}\right) \mu'_1(d_1)
\mu_2\left(\frac{-\det(k_2)}{\pi^{-s}c_2}\right) \left(\frac{\alpha_2}{\beta_2}\right)^s
& \text{, if }\mu_1
\text{ is ramified } \text{ and }   \mu'_2 \text{ is unramified ; } \\ 
\mu'_1(d_1)
\mu_2\left(\frac{-\det(k_2)}{\pi^{-s}c_2}\right) \left(\frac{\alpha_2}{\beta_2}\right)^s
& \text{, if } \mu_1
\text{ and }   \mu'_2 \text{ are unramified. }
\end{cases}$$
\end {lemma}

\noindent{\it Proof : } 
By definition $F(k_1,k_2)=0$ unless there exist
$k_0=\begin{pmatrix} 1 & b_0 \\ c_0 & 1  \end{pmatrix}
\in I_{f}$  such that 
$$k_1k_0^{-1}\in B \qquad \text{ and } \qquad k_2k_0^{-1}\begin{pmatrix} 0 & 1 \\ 1 & 0
\end{pmatrix}\in B ,$$ in which case
$$F(k_1,k_2)=\chi_1(k_1k_0^{-1})\chi_2 \Bigl( k_2k_0^{-1}\begin{pmatrix} 0 & 1 \\ 1 & 0 \end{pmatrix} \Bigr) 
\delta^{\frac{1}{2}}\Bigl(k_1k_0^{-1}k_2k_0^{-1}\begin{pmatrix} 0 & 1 \\ 1 & 0 \end{pmatrix} \Bigr) f(k_0).$$
From $k_1k_0^{-1}\in B$, we deduce that $c_1=c_0d_1$. 
From $k_2k_0^{-1}\begin{pmatrix} 0 & 1 \\ 1 & 0 \end{pmatrix}\in B $  
we deduce that $d_2=b_0c_2$. Hence 
$$d_1\in \OF^\times, \qquad  \frac{c_1}{d_1}\in \pi^x{\OF}^{\mu_1}, \qquad c_2 \neq 0 
\qquad \text{and} \qquad \frac{d_2}{c_2}\in \pi^{-y}{\OF}^{\mu'_2}.$$ 
Moreover  
$$k_1k_0^{-1}=\begin{pmatrix}\frac{\det{k_1}}{d_1\det{k_0}} & * \\ 0 & d_1 \end{pmatrix}
\text{ and } 
k_2k_0^{-1}\begin{pmatrix} 0 & 1 \\ 1 & 0\end{pmatrix}
=\begin{pmatrix}  \frac{-\det{k_2}}{c_2\det{k_0}} & * \\ 0 & c_2 \end{pmatrix}.$$
Since $x-y\geq n_1-m_1$, $x-y \geq n_2-m_2$ and  $x-y \geq 1$ we have
$$\mu_1(\det{k_0})=\mu_2(\det{k_0})=1.$$ 
Hence
$$
F(k_1,k_2) =  \mu_1(\frac{\det{k_1}}{d_1})
\mu'_1(d_1)\mu_2(\frac{-\det{k_2}}{c_2}) \mu'_2(c_2)
\left\vert \frac{1}{c_2}\right\vert f\begin{pmatrix} 1 & \frac{d_2}{c_2}\\
\frac{c_1}{d_1}  & 1 \end{pmatrix}.$$

From here and (\ref{f}) follows the desired formula for $F$. 

\medskip
Conversely, if  $k_1$ and $k_2$ are such that 
 $\frac{c_1}{d_1}\in \pi^x{\OF}^{\mu_1}$ and 
 $\frac{d_2}{c_2}\in \pi^{-y}{\OF}^{\mu'_2}$ one can take 
$$k_0=\begin{pmatrix} 1  & d_2c_2^{-1} \\ c_1d_1^{-1} & 1  \end{pmatrix}.$$
\hfill $\Box$

\begin{rem}
One can compute $F$ without the assumption $x-y \geq\max (n_1-m_1, n_2-m_2,1)$.
However, $F$ needs not decompose as 
a product of functions of one variable as in the above lemma. 

For example, if $x=n_3=0$ and  $n_1=n_2$, then for all $k_1\in K$ and  $k_2\in K$ 
$$F(k_1,k_2)=\begin{cases}
   \omega_1(\frac{c_1d_2-d_1c_2}{\det k_2})& \text{ , if } d_1\in \OF^\times, 
  \enspace c_2\in \OF^\times \text{ and } c_1d_2\neq d_1c_2 \\
0 & \text{ ,  otherwise}. \end{cases}$$
\end{rem}

\subsection{From $\Psi$ to $\ell$}

Now, we want to express $F \in V_1\otimes V_2$
in terms of the new vectors $v_1$ and $v_2$. 

  From now on we suppose that $x$, $y$ and $z$ are integers as in theorem \ref{no-vt}.
  We may also suppose that $x\geq 1$, because otherwise $V_1$, $V_2$ and $V_3$
are all unramified and this case is covered in Theorem \ref{vt-000}.
Observe also that if $y=0$, then $\mu'_2$ is unramified and therefore 
$ \OF^{\mu'_2}=\OF$.

For $i=1,2$, since $k_i\in K$, both $c_i$ and $d_i$ are in $\OF$, 
and one of them is in $\OF^\times$. Hence 
\begin{itemize}
\item $\frac{c_1}{d_1}\in \pi^x{\OF}^\times$ if, and only if  $k_1\in  I_{x}\backslash I_{x+1}$, 
\item $\frac{c_1}{d_1}\in \pi^x{\OF}$ if, and only if  $k_1\in  I_{x}$, 
\item $\frac{d_2}{c_2}\in \pi^{-y}{\OF}^\times$ with $y\geq 1$ if, and only if  $k_2\in I_{y}\backslash I_{y+1}$, 
\item $\frac{d_2}{c_2}\in \pi^{-y}{\OF}$ with $y\geq 0$ if, and only if  $k_2\in K \backslash I_{y+1}$.
\end{itemize} 

\begin{lemma} \label{calcul-F}
With the notations of (\ref{v12}), $F$ is a non-zero multiple of $ v_1^* \otimes v_2^*$.
\end{lemma}

{\it Proof : } Both $F$ and $ v_1^* \otimes v_2^*$ are elements in 
$\Ind_{B \times B}^{G \times G} \Bigl( \chi_1 \times \chi_2
\Bigr)$, hence it is enough to compare their restrictions to $K\times
K$. By  the above discussion together with lemmas \ref{FV} and \ref{support} the two restrictions are supported by 
$$\begin{cases} 
(I_{x}\backslash I_{x+1})\times (I_{y}\backslash I_{y+1}) &  
\text{ , if } \mu_1 \text{ and }   \mu'_2 \text{ are ramified ; }\\
I_{x}\times (I_{y}\backslash I_{y+1}) &  
\text{ , if } \mu_1 \text{ is unramified } \text{ and }   \mu'_2 \text{ is ramified ;  }\\
(I_{x}\backslash I_{x+1})\times (K \backslash I_{y+1}) &  
\text{ , if } \mu_1 \text{ is ramified } \text{ and }   \mu'_2 \text{ is unramified ; }\\
I_{x}\times  (K \backslash I_{y+1}) & 
\text{ , if } \mu_1 \text{ and }   \mu'_2 \text{ are unramified. }
\end{cases}$$

There are $16$  different cases depending on whether each one 
among  $\mu_1$,  $\mu'_1$,   $\mu_2$ and  $\mu'_2$ is ramified or
unramified. Since it is a straightforward verification from lemmas \ref{calcul-NR}, \ref{calcul-SP0}, \ref{calcul-SP1} and \ref{calcul-SP2}, in order 
to avoid repetitions or cumbersome notations, we will only give
the final result : 
\begin{equation}\begin{split}
 &F=\lambda_1\lambda_2\mu_2(-1) \alpha_1^{m_1-x}\alpha_2^{m_2}\beta_2^{-y}
 (v_1^* \otimes v_2^*) \text{ , where }\\
&\lambda_i=\begin{cases} \Bigl(1-\frac{\beta_i}{\alpha_i}\Bigr)^{-1} &  \text{ , if }  V_i \text{ is unramified,}\\
1  & \text{ , if }  V_i \text{ is ramified. }
\end{cases}  
 \end{split}\end{equation}

In all cases $F$ is a non-zero multiple of $ v_1^* \otimes v_2^*$.
\hfill $\square$

\bigskip
Since by definition $\ell(F\otimes\bullet)=\Psi(F)=\Phi(f)$, the above 
lemma together with  lemma \ref{phi-f-v_3} imply  theorem \ref{no-vt}.

\subsection{Proof of Theorems  \ref{vt-00n} and \ref{vt-mkn}.}

We assume henceforth that  $\mu_1$ and $\mu'_2$ are {\it both}
unramified ($n_1-m_1=m_2=0$). We put $y=z=0$ and
 $x=\max(n_1,n_3)\geq 1$. Since $\omega_1\omega_2\omega_3=1$,  
$\max(n_1,n_3)=\max(n_1,n_2,n_3)\geq 1$.

Then lemma \ref{phi-f-v_3}  yields :
\begin{equation}
\ell(F\otimes v_3)=\Psi(F)(v_3)=\Phi(f)(v_3)\neq 0.  
\end{equation}

From this and lemma \ref{calcul-F} we deduce :
\begin{lemma}\label{l-F-v_3} 
We have $\ell(v_1^*\otimes v_2^*\otimes v_3)\neq 0$ where
\begin{equation*}\begin{split}
&v_1^*=\begin{cases} \gam^{x-n_1}\cdt v_1 &  \text{ , if }  \mu'_1 \text{ is ramified,}\\
\gam^{x}\cdt v_1-\beta_1 \gam^{x-1}\cdt v_1  & \text{ , if }  \mu'_1 \text{ is unramified.}
\end{cases}  \\
&v_2^*=\begin{cases}  v_2 &  \text{ , if }  \mu_2 \text{ is ramified.}\\
v_2-\alpha_2^{-1} \gam\cdt v_2  & \text{ , if }  \mu_2 \text{ is unramified.}
\end{cases} \end{split}\end{equation*}
\end{lemma}

\subsubsection{The case of two unramified representations.}

Suppose that $n_1=n_2=0$, so that  $x=n_3$. Then   lemma \ref{l-F-v_3} yields : 
$$\ell\Bigl((\gam^{{n_3}}\cdt v_1-\beta_1\gam^{n_3-1}\cdt v_1)\otimes
(\gam\cdt v_2-\alpha_2v_2)\otimes v_3\Bigr)\neq 0. $$

This expression can be simplified as follows. Consider for $m\geq 0$ 
the linear form :
$$\psi_m(\bullet)=\ell(\gam^{m}\cdt v_1\otimes  v_2 \otimes\bullet)\in\widetilde{V_3}.$$ 

As observed in the introduction, 
$\psi_m$ is invariant by $\gam^{m}K\gam^{-m}\cap K =I_m $, hence vanishes 
if $m<n_3=\cond(\widetilde{V_3})$. Therefore, for ${n_3}\geq 2$ : 
\begin{equation*}\begin{matrix}
\ell\Bigl((\gam^{{n_3}}\cdt v_1-\beta_1\gam^{n_3-1}\cdt v_1)\otimes
(\gam\cdt v_2-\alpha_2v_2)\otimes v_3\Bigr)\hfill \\
\hskip1cm = -\alpha_2\psi_{n_3}(v_3)+\beta_1\alpha_2\psi_{{n_3}-1}(v_3)+\psi_{{n_3}-1}(\gam^{-1} \cdt v_3)-\beta_1\psi_{{n_3}-2}(\gam^{-1} \cdt v_3)\hfill\\
\hskip1cm =  -\alpha_2\psi_{n_3}(v_3)\hfill\\
\hskip1cm =  -\alpha_2\ell(\gam^{{n_3}}\cdt v_1\otimes  v_2\otimes v_3) \neq 0.\hfill\\
\end{matrix}\end{equation*}

If ${n_3} = 1$, only the two terms in the middle vanish and we obtain
$$\alpha_2\ell(\gam\cdt v_1\otimes  v_2\otimes v_3)
+\beta_1\ell(v_1\otimes \gam\cdt v_2\otimes v_3)\neq 0.$$

Put 
$g=\begin{pmatrix}0 & 1 \cr \pi & 0\cr\end{pmatrix}$. Then 
$g \gam = \begin{pmatrix} 0 & 1 \\ 1 & 0  \end{pmatrix}\in K$ 
and $\gam^{-1} g = \begin{pmatrix} 0 & \pi \\ \pi & 0  \end{pmatrix}\in \pi
K$. Hence : 
$$\begin{matrix}
\beta_1\ell(v_1\otimes \gam\cdt v_2\otimes v_3)\hfill 
& = & \beta_1\ell(\gam\gam^{-1} g\cdt v_1\otimes g\gam\cdt v_2\otimes g\cdt v_3)\hfill \\
& = & \beta_1\omega_1(\pi)\ell\bigl(\gam\cdt v_1\otimes  v_2\otimes g\cdt v_3\bigr)\hfill \\
& = & \alpha_1^{-1}\ell(\gam\cdt v_1\otimes  v_2\otimes g\cdt v_3).\hfill \\
\end{matrix}$$

Therefore 
$$\ell\Bigl(\gam\cdt v_1\otimes v_2\otimes(g\cdt v_3+\alpha_1\alpha_2v_3)\Bigr)\neq 0,$$ 
in particular $$\Psi(\gam\cdt v_1\otimes v_2)\neq 0.$$
By the same argument as in paragraph \ref{cas-simple2}  we conclude 
that 
$$\ell(\gam\cdt v_1\otimes v_2\otimes v_3)=\Psi(\gam\cdt v_1\otimes v_2)(v_3)\neq 0.$$ 
Hence, if  $n_3\geq 1$, $\gam^{{n_3}}\cdt v_1\otimes v_2\otimes v_3$
is a test vector. This completes the proof of Theorem \ref{vt-00n}.

\subsubsection{The case of two ramified principal series.}

Suppose that $V_1$ and $V_2$ are both ramified ($m_1>0$, $n_1-m_1=0$, $m_2=0$, $n_2>0$) and put $n=x-n_1=\max(n_2-n_1,n_3-n_1)$. Then   lemma \ref{l-F-v_3} yields : 
$$\ell(\gam^{n}\cdt v_1 \otimes v_2   \otimes v_3)\neq 0,$$ 
hence $\gam^{n}\cdt v_1 \otimes v_2   \otimes v_3$ is a test vector.

\subsubsection{The case where $V_1$ is unramified and  $V_2$ is  ramified.}

Suppose that $n_1=0$, but  $n_2>0$. Then  $x=n_3\geq n_2$ and  
lemma \ref{l-F-v_3} yields : 

$$\ell\Bigl((\gam^{{n_3}}\cdt v_1-\beta_1\gam^{n_3-1}\cdt v_1)\otimes v_2\otimes v_3\Bigr)\neq 0.$$

If $n_2<{n_3}$, then $$\gam^{n_3-1}K\gam^{1-n_3}\cap I_{n_2}\supset I_{{n_3}-1},$$
 and therefore $$\ell(\gam^{n_3-1}\cdt v_1\otimes v_2\otimes
\bullet)\in \widetilde{V_3}^{I_{{n_3}-1},\omega_3^{-1}}=\{0\}.$$ 
Hence 
$$\ell(\gam^{{n_3}}\cdt v_1\otimes v_2\otimes v_3)\neq 0,$$ 
that is $\gam^{{n_3}}\cdt v_1\otimes v_2\otimes v_3$ is a test vector.

If $n_2={n_3}$, the condition on the central character 
forces $V_3$ and $\omega_3$ to have the same conductor. Hence $V_3$ is also a principal series. 
In this case we do not see {\it a priori} a reason for either 
$\ell(\gam^{{n_3}}\cdt v_1\otimes v_2\otimes v_3)$ or 
$\ell(\gam^{n_3-1}\cdt v_1\otimes v_2\otimes v_3)$ to vanish. But we can notice that the two linear forms 
$$\ell(\gam^{{n_3}}\cdt v_1\otimes v_2\otimes \bullet)  \quad {\rm and} \quad \ell(\gam^{n_3-1}\cdt v_1\otimes v_2\otimes \bullet) $$
belong both to the new line $\widetilde{V_3}^{I_n,\omega_3^{-1}}$ of $\widetilde{V_3}$, hence they are proportionals.

\subsubsection{The case where $V_1$ is ramified and  $V_2$ is  unramified.}

Suppose that $n_1>0$  and $n_2=0$. Then  $x=n_3\geq n_1$ and  
lemma \ref{l-F-v_3} yields : 
$$\ell\Bigl(\gam^{n_3-n_1}\cdt v_1 \otimes (\gam\cdt v_2-\alpha_2 v_2)\otimes v_3\Bigr)\neq 0. $$
 
If $n_1<{n_3}$, then 
$$\ell(\gam^{n_3-n_1-1}\cdt v_1\otimes v_2\otimes \bullet)\in \widetilde{V_3}^{I_{{n_3}-1},\omega_3^{-1}}=\{0\}.$$ 
Then 
$$\ell(\gam^{n_3-n_1}\cdt v_1 \otimes \gam\cdt v_2\otimes v_3)=
\ell(\gam^{n_3-n_1-1}\cdt v_1 \otimes v_2\otimes\gam^{-1} \cdt v_3)=0.$$
Hence 
$$\ell(\gam^{n_3-n_1}\cdt v_1 \otimes v_2\otimes v_3)\neq 0,$$ 
that is $\gam^{n_3-n_1}\cdt v_1 \otimes v_2\otimes v_3$ is a test vector.

If $n_1={n_3}$, the condition on the central character 
forces $V_3$ to be also a principal series. In this case 
we do not see {\it a priori} a reason for either 
$\ell(v_1 \otimes v_2\otimes v_3)$ or 
$\ell(v_1 \otimes \gam\cdt v_2\otimes v_3)$ to vanish. 
But we can once again notice that the two linear forms
$$\ell(v_1 \otimes v_2\otimes \bullet) \quad {\rm and} \quad  \ell(v_1
\otimes  \gam\cdt v_2\otimes \bullet)$$
belong to the line generated by a new vector in $\widetilde{V_3}$, 
hence are proportionals.

The proof of Theorem \ref{vt-mkn} is now complete.

\end{document}